\documentclass[a4paper,reqno]{amsart}

\usepackage{graphicx}

\usepackage{amsthm}
\theoremstyle{plain}
\newtheorem{theorem}{Theorem}

\newtheorem{corollary}[theorem]{Corollary}
\theoremstyle{definition}

\newtheorem{example}{Example}
\theoremstyle{remark}
\newtheorem{remark}[theorem]{Remark}

\newcommand{\emphbf}[1]{\textbf{\textit{#1}}}

\newcommand{\R}{\mathbb{R}} 
\newcommand{\F}{\mathbb{F}}
\newcommand{\T}{\mathbb{T}}

\newcommand{\Ad}{\operatorname{Ad}}

\newcommand{\G}{\mathfrak{g}}

\begin{document}
\bibliographystyle{plain}

\title[Momentum and energy nonholonomic integrators]{Momentum and energy preserving
integrators for nonholonomic dynamics}

\author[S.\ Ferraro]{S.\ Ferraro}
\address{S.\ Ferraro:
Instituto de Ciencias Matem\'aticas (CSIC-UAM-UC3M-UCM), Serrano 123, 28006
Madrid, Spain} \email{sferraro@uns.edu.ar}

\author[D.\ Iglesias]{D.\ Iglesias}
\address{D.\ Iglesias:
Instituto de Ciencias Matem\'aticas (CSIC-UAM-UC3M-UCM), Serrano 123, 28006
Madrid, Spain} \email{iglesias@imaff.cfmac.csic.es}

\author[D.\ Mart\'{\i}n de Diego]{D.\ Mart\'{\i}n de Diego}
\address{D.\ Mart\'{\i}n de Diego:
Instituto de Ciencias Matem\'aticas (CSIC-UAM-UC3M-UCM), Serrano 123, 28006
Madrid, Spain} \email{d.martin@imaff.cfmac.csic.es}

\begin{abstract}
In this paper,  we propose a geometric integrator for nonholonomic
mechanical systems. It can be applied to discrete Lagrangian
systems specified through a discrete Lagrangian $L_d\colon Q\times
Q\to \mathbb{R}$, where $Q$ is the configuration manifold, and a
(generally nonintegrable) distribution $\mathcal{D}\subset TQ$. In
the proposed method, a discretization of the constraints is not
required. We show that the method preserves the discrete
nonholonomic momentum map, and also that the nonholonomic
constraints are preserved in average. We study in particular the
case where $Q$ has a Lie group structure and the discrete Lagrangian and/or
nonholonomic constraints have various invariance properties, and
show that the method is also energy-preserving in some important
cases.
\end{abstract}

\maketitle

\tableofcontents

\section{Introduction}

During the last years, there has been an increasing interest in
nonholonomic mechanical systems, in part motivated by some open
questions in the subject, such as those concerning reduction,
integrability, stabilization or controllability; and also for
their applicability in engineering, specially in robotics, mainly
since it describes the motion of wheeled devices (see
\cite{Bloch:Nonholonomic_mechanics_and_control,
Bloch_Krishnaprasad_Marsden_Ratiu:The_Euler-Poincare_equations_and_double_bracket_dissipation,
de_Leon_Martin_de_Diego:On_the_geometry_of_non-holonomic_Lagrangian_systems}
and the expository paper
\cite{Bloch_Marsden_Zenkov:Nonholonomic_dynamics}).

When a mechanical system is subjected to some external
constraints, the latter may be expressed in terms of relations
imposing restrictions on the allowable positions and velocities.
The constraints are then called \emph{nonholonomic} if the velocity
dependence is essential, in the sense that the constraint
relations can not be reduced, by integration, to relations
depending on the position coordinates only. Geometrically,
nonholonomic constraints are globally described by a submanifold
${\mathcal D}$ of the velocity phase space $TQ$. In most of
the known examples ${\mathcal D}$ is a vector subbundle of $TQ$,
i.e., the constraints have a linear dependence on the velocities.
 Lagrange--d'Alembert's
principle allow us to determine the set of possible values of the
constraint forces from the constraint manifold ${\mathcal D}$.
Then, to determine the dynamics of the nonholonomic system, it is only necessary  to fix
initially the pair $(L, {\mathcal D})$, where $L\colon TQ\to
\R$ is  a Lagrangian function, usually of mechanical type (see
\cite{Bloch:Nonholonomic_mechanics_and_control,
Bullo_Lewis:Geometric_control_of_mechanical_systems,
Cendra_Ibort_de_Leon_de_Diego:A_generalization_of_Chetaevs_principle_for_a_class_of_higher_order_nonholonomic_constraints}
for an extension of the classical Lagrange--d'Alembert's
principle).

Very recently, many authors
\cite{Cortes_Martinez:Non-holonomic_integrators,
de_Leon_Martin_de_Diego_Santamaria-Merino:Geometric_numerical_integration_of_nonholonomic_systems_and_optimal_control_problems,
Fedorov_Zenkov:Discrete_nonholonomic_LL_systems_on_Lie_groups,
Iglesias_Marrero_MartindeDiego_Martinez:Discrete_Nonholonomic_Lagrangian_Systems_on_Lie_Groupoids,
McLachlan_Perlmutter:Integrators_for_nonholonomic_mechanical_systems}
started the study of geometric integrators adapted to nonholonomic
systems, obtaining very stable numerical integrators with some
preservation properties (such as discrete nonholonomic momentum
map preservation) and very good energy behavior. This problem is
of
 considerable interest given the crucial role of nonholonomic dynamics in many applications in engineering.
{}From the numerical point of view, in
\cite{McLachlan_Scovel:A_survey_of_open_problems_in_symplectic_integration}
it appeared as an open question:  ``...The problem for the more
general class of non-holonomic constraints is still open, as is
the question of the correct analogue of symplectic integration for
non-holonomically constrained Lagrangian systems...''.

 The most
interesting approach to nonholonomic integrators appears as an
adaptation of the so-called variational integrators
\cite{Marsden_West:Discrete_mechanics_and_variational_integrators}
incorporating  a discrete constraint submanifold, in addition to
a discretization of the Lagrangian function and the vector
subbundle $\mathcal{D}$. Then, the numerical method is obtained
from the so-called Discrete Lagrange--d'Alembert's principle~\cite{Cortes_Martinez:Non-holonomic_integrators}, recovering many
of the geometric properties of the continuous system.

Obviously, since nonholonomic mechanics is not symplectic-preserving, it seems interesting to try to preserve another
geometric invariance property of the continuous nonholonomic
system, as for instance, the energy function in the autonomous
case. This is precisely the starting point of view of our paper.
Moreover, a discretization of the constraints is not required
here. We show that the method preserves the discrete nonholonomic
momentum map, and also that the constraints are preserved in
average. We study in particular the case where the configuration
space is a Lie group and the discrete Lagrangian and/or
nonholonomic constraints have various invariance properties, and
show that the method is also energy-preserving in many important
cases. In particular, the main result of the  paper, Theorem~\ref{thm:preservation of energy}, states that
if the configuration
space is a Lie group and the Lagrangian is defined by a bi-invariant Riemannian metric, then, from a
left-invariant discretization of the Lagrangian, we obtain a
\emphbf{fixed time-step, energy-preserving numerical method} for the
continuous nonholonomic system, without requiring any invariance conditions on $\mathcal{D}$.
See~\cite{Cortes:Energy_conserving_nonholonomic_integrators} for a
variable time-step algorithm that preserves energy.

The paper is structured as follows. In Section \ref{section2}, we
introduce continuous nonholonomic mechanical systems for the
case of mechanical energy Lagrangians defined by a given
Riemannian metric and a potential function. In this case, the
equations of motion for the constrained system are geodesic
equations for an affine connection (in the kinetic case) that is
not generally Levi-Civita, obtained from the induced orthogonal
projection onto the nonholonomic distribution
 (see
\cite{Cantrijn_Cortes_deLeon_MartindeDiego:On_the_geometry_of_generalized_Chaplygin_systems,
Lewis:Affine_connections_and_distributions_with_applications_to_nonholonomic_mechanics}).
 In Section \ref{Dvc} we recall some  definitions concerning discrete variational mechanics
 (discrete Lagrangian, discrete Euler--Lagrange equations, discrete flow, momentum map...).
The new proposed method appears in Section \ref{section3},
constructed from the discrete Lagrangian and the orthogonal
projectors induced by the distribution $\mathcal{D}$ and the Riemannian metric. Then we consider the
case when the configuration space is a Lie group and we obtain
under adequate invariance properties the preservation of energy. In addition, we study the momentum nonholonomic map for the
proposed nonholonomic integrator. In Section \ref{section5}, we introduce a nonholonomic version
of the St\"ormer--Verlet method which is a natural extension of the RATTLE method for nonholonomic systems.
In  Section  \ref{section6} we test our
method in three examples (the nonholonomic particle, the
snakeboard and the Chaplygin sleigh). 
The paper ends with a section of
conclusions and future work.

\section{Continuous nonholonomic mechanics}\label{section2}
We shall start with a configuration space $Q$, which is an
$n$-dimensional differentiable manifold with local coordinates
$(q^i)$, $1\leq i\leq n=\dim Q$. Constraints linear in the
velocities are given by equations of the form
\[
\phi^{a}(q^i, \dot{q}^i)=\mu^a_i(q)\dot{q}^i=0, \quad 1\leq a\leq
m,
\]
depending, in general, on configuration coordinates and  their velocities.
{}From an intrinsic point of view, the linear
constraints are defined by a distribution
${\mathcal D}$ on $Q$ of rank $n-m$ such that the
annihilator of ${\mathcal D}$ is locally given by
\[
{\mathcal D}^o = \operatorname{span}\left\{ \mu^{a}=\mu_i^{a}dq^i \; ; 1 \leq a
\leq m \right\}
\]
where the one-forms $\mu^{a}$ are independent.

The various kinds of constraints we are concerned with will roughly come in two types: \emphbf{holonomic} and \emphbf{nonholonomic},
depending on whether  the constraint  is derived from a constraint in the configuration space or not.
Therefore, the dimension of the space of configurations is reduced by  holonomic constraints but not by nonholonomic constraints.
Thus, holonomic constraints allow a reduction in the number of coordinates of the configuration space needed to formulate a given problem
(see~\cite{Neimark_Fufaev:Dynamics_of_Nonholonomic_Systems}).

We will restrict ourselves to the case  of nonholonomic
constraints. In this case, the constraints are given by a
nonintegrable distribution $\mathcal{D}$. In addition to these constraints,
we need to specify the dynamical evolution of the system, usually by
fixing a Lagrangian function $L\colon  TQ \to \R$. In mechanics,
the central concepts permitting the extension of mechanics from
the Newtonian point of view to the Lagrangian one are the notions
of virtual displacements and virtual work; these concepts were
formulated in the developments of mechanics, in their application
to statics. In nonholonomic dynamics,  the procedure is given by
the \emphbf{Lagrange--d'Alembert principle}.
 This principle allows us to determine the set of possible values of the constraint forces from the set $\mathcal{D}$ of admissible kinematic states alone. The resulting equations of motion are
\begin{equation*}
\left[ \frac{d}{dt}\left( \frac{\partial L}{\partial \dot
q^i}\right) - \frac{\partial L}{\partial q^i} \right] \delta
q^i=0,
\end{equation*}
where $\delta q^i$ denotes the virtual displacements verifying
\begin{equation*}
\mu^a_i\delta q^i =0
\end{equation*}
(for the sake of simplicity, we will assume that the system is not subject to non-conservative forces). This must be supplemented by the constraint equations.
By using the Lagrange multiplier rule, we obtain
\begin{equation*}
\frac{d}{dt}\left( \frac{\partial L}{\partial \dot
q^i}\right)-\frac{\partial L}{\partial q^i}={\lambda}_a\mu^a_i  .
\end{equation*}
The term on the right represents the  constraint force or reaction force induced by the constraints.
The functions $\lambda_a$ are Lagrange multipliers which, after being computed using the constraint equations, allow us to obtain a set of second order differential equations.

Now we restrict ourselves to the case of  nonholonomic mechanical
systems where the Lagrangian is of mechanical type
\[
L(v_q)=\frac{1}{2}g(v_q, v_q) - V(q), \quad v_q\in T_qQ,
\]
where $g$ is a Riemannian metric on the configuration space $Q$.
Locally, the metric is determined by the matrix $M=(g_{ij})_{1\leq
i, j\leq n}$ where $g_{ij}=g(\partial/\partial q^i,
\partial/\partial q^j)$.

Using some basic tools of Riemannian geometry, we may write the  equations of motion of the unconstrained
system as
\begin{equation}\label{free}
\nabla _{\dot{c}(t)} \dot{c}(t) = - {\rm grad}~ V(c(t)),
\end{equation}
where $\nabla$ is the Levi--Civita connection associated to $g$.
Observe that if $V\equiv 0$ then the Euler--Lagrangian equations
are the equations of the geodesics for the Levi-Civita connection.

When the system is subjected to nonholonomic constraints, the equations become
\begin{equation*}
\nabla _{\dot{c}(t)} \dot{c}(t) = - {\rm grad}~
V(c(t))+\lambda(t),\quad \dot{c}(t) \in
{\mathcal D}_{c(t)},
\end{equation*}
where  $\lambda$ is a section of ${\mathcal D}^{\perp}$ along $c$.
Here ${\mathcal D}^{\perp}$ stands for the orthogonal complement
of ${\mathcal D}$ with respect to the metric $g$.

In coordinates, by defining the $n^3$ functions $\Gamma^k_{ij}$ (Christoffel symbols for $\nabla$) by
\[
\nabla_{\!\!\frac{\partial}{\partial
q^i}}\,\frac{\partial}{\partial
q^j}=\Gamma^k_{ij}\frac{\partial}{\partial q^k},
\]
we may rewrite the nonholonomic equations of motion as
\begin{align*}
\ddot{q}^k(t)+\Gamma^k_{ij}(c(t))\dot{q}^i(t)\dot{q}^j(t)&=-g^{ki}(c(t))\frac{\partial V}{\partial q^i}+ \lambda_a(t)g^{ki}(c(t))\mu^a_i(c(t)) 
\\
\mu^a_i(c(t))\dot{q}^i(t)&=0 
\end{align*}
where $t\mapsto (q^1(t), \ldots, q^n(t))$ is the local
representative of $c$ and $(g^{ij})$ is the inverse matrix of $M$.

Since $g$ is a Riemannian metric, the $m \times
m$ matrix $(C^{ab})=(\mu^a_ig^{ij}\mu^b_j)$ is
symmetric and regular.
 Define now the vector fields $Z^a$, $1\leq a\leq m$ on $Q$
 by
\[
g(Z^{a}, Y)=\mu^a(Y), \ \hbox{ for all vector
fields } Y, \; 1 \leq a \leq m;
\]
that is, $Z^{a}$ is the gradient vector field of the 1-form
$\mu^{a}$. Thus, ${\mathcal D}^{\perp}$ is spanned by $Z^{a}$, $1\leq a\leq m$. In local
coordinates, we have
\[
Z^{a} = g^{ij} \mu^{a}_i \frac{\partial}{\partial q^j}.
\]

We can construct two complementary
projectors
\begin{align*}
{\mathcal P}\colon &TQ\to {\mathcal D}\\
{\mathcal Q}\colon &TQ\to {\mathcal
D}^{\perp},
\end{align*}
orthogonal with respect to the metric $g$. The projector ${\mathcal Q}$
is locally described by
\[
{\mathcal Q} = C_{ab} Z^a \otimes \mu^b= C_{ab}g^{ij} \mu^{a}_i
\mu^b_k\frac{\partial}{\partial q^j} \otimes dq^k.
\]
Using these projectors we may rewrite the equations of motion as
follows.
A curve $c(t)$ is a motion for the nonholonomic system if it
satisfies the constraints, i.e., $\dot{c}(t)\in {\mathcal
D}_{c(t)}$, and, in addition, the ``projected equation of motion''
\begin{equation}\label{q5}
{\mathcal P}(\nabla_{\!\dot{c}(t)} \dot{c}(t)) =
- {\mathcal P}(\hbox{grad} ~ V(c(t)))
\end{equation}
is fulfilled.

Summarizing, we have obtained the dynamics of the nonholonomic
system~\eqref{q5} applying the projector ${\mathcal P}$ to the
dynamics of the free system~\eqref{free}. In Section~\ref{section3}, we will
 use $\mathcal{P}$ and $\mathcal{Q}$ to obtain a geometric integrator for nonholonomic systems.

\section{Variational integrators}
\label{Dvc}

The equations of motion for an unconstrained Lagrangian system given by a
Lagrangian function $L\colon TQ\to \R$ are the well-known
Euler--Lagrange equations \begin{equation*}
\frac{d}{dt}\left(\frac{\partial L}{\partial\dot
q^i}\right)-\frac{\partial L}{\partial q^i}=0, \quad 1\leq i\leq
n.
\end{equation*}
 It is well known that the origin of these equations is variational (see \cite{Abraham_Marsden:Foundations_of_Mechanics}).
 Now, variational integrators retain  this variational character
and also some of the geometric properties of the continuous
system, such as symplecticity and momentum conservation (see
\cite{Hairer_Lubich_Wanner:Geometric_numerical_integration,Marsden_West:Discrete_mechanics_and_variational_integrators}
and references therein).

In the following we will summarize the main features of this type
of numerical integrators.  A \emphbf{discrete Lagrangian} is a map
$L_d\colon Q \times Q\to  \R$, which may be considered as
an approximation of a continuous Lagrangian $L\colon TQ\to
\R$. Define the \emphbf{action sum} $S_d\colon Q^{N+1}\to
\R$ corresponding to the Lagrangian $L_d$ by
\[
{S_d}=\sum_{k=1}^{N}  L_d(q_{k-1}, q_{k}),
\]
where $q_k\in Q$ for $0\leq k\leq N$. The discrete variational
principle   states that the solutions of the discrete system
determined by $L_d$ must extremize the action sum given fixed
endpoints $q_0$ and $q_N$. By extremizing ${S_d}$ over $q_k$,
$1\leq k\leq N-1$, we obtain the system of difference equations
\begin{equation}\label{discreteeq}
 D_1L_d( q_k, q_{k+1})+D_2L_d( q_{k-1}, q_{k})=0.
\end{equation}
or, in coordinates,
\[
\frac{\partial L_d}{\partial q_0^i}(q_k, q_{k+1})+\frac{\partial
L_d}{\partial q_1^i}(q_{k-1}, q_{k})=0, \ 1\leq i\leq n,\ 1\leq
k\leq N-1.
\]

These  equations are usually called the  \emphbf{discrete
Euler--Lagrange equations}. Under some regularity hypotheses (the
matrix $(D_{12}L_d(q_k, q_{k+1}))$ is regular), it is possible to
define a (local) discrete flow $ \Upsilon\colon Q\times
Q\to  Q\times Q$, by $\Upsilon(q_{k-1}, q_k)=(q_k,
q_{k+1})$ from~\eqref{discreteeq}. Define the  discrete
Legendre transformations associated to  $L_d$ as
\begin{align*}
\F^-L_d\colon Q\times Q&\to  T^*Q\\
(q_0, q_1)&\longmapsto (q_0, -D_1 L_d(q_0, q_1))\\
\F^+L_d\colon  Q\times Q&\to  T^*Q\\
(q_0, q_1)&\longmapsto (q_1, D_2 L_d(q_0, q_1))\; ,
\end{align*}
and the discrete Poincar{\'e}--Cartan 2-form $\omega_d=(\F^+L_d)^*\omega_Q=(\F^{-}L_d)^*\omega_Q$,
where $\omega_Q$ is the canonical symplectic form on $T^*Q$. The
discrete algorithm determined by $\Upsilon$ preserves the
symplectic form $\omega_d$, i.e., $\Upsilon^*\omega_d=\omega_d$.
Moreover, if the discrete Lagrangian is invariant under the
diagonal action of a Lie group $G$, then the discrete momentum map
$J_d\colon Q\times Q \to  {\mathfrak g}^*$ defined by
\[ \langle
J_d(q_k, q_{k+1}), \xi\rangle=\langle D_2L_d(q_k, q_{k+1}),
\xi_Q(q_{k+1})\rangle \]
is preserved by the discrete flow.
Therefore, these integrators are symplectic-momentum preserving. Here, $\xi_Q$ denotes the fundamental vector field
determined by $\xi\in {\mathfrak g}$, where ${\mathfrak g}$ is the Lie
algebra of $G$.

\section{A geometric nonholonomic integrator}\label{section3}
This work proposes a numerical method for the integration of nonholonomic systems. It is not truly variational; however, it is geometric in nature and we show in Corollary~\ref{cor:preservation of momentum} that it preserves the discrete nonholonomic momentum map in the presence of horizontal symmetries. Moreover, we prove in Theorem~\ref{thm:preservation of energy} that under certain symmetry conditions, the energy of the system is preserved.

Consider a discrete Lagrangian $L_d\colon Q\times Q\to \R$. The proposed discrete nonholonomic equations are
\begin{subequations}\label{eq:propuesta original}
\begin{align}
{\mathcal P}^*_{q_k}( D_1 L_d(q_k, q_{k+1}))+
{\mathcal P}^*_{q_k} (D_2 L_d (q_{k-1},
q_k))&=0 \\
{\mathcal Q}^*_{q_k}( D_1 L_d(q_k, q_{k+1}))-
{\mathcal Q}^*_{q_k} (D_2 L_d (q_{k-1},
q_k))&=0,
\end{align}
\end{subequations}
where the subscript $q_k$ emphasizes the fact that the projections take place in the fiber over $q_k$.
The first equation is the projection of the discrete
Euler--Lagrange equations to the constraint distribution
$\mathcal{D}$, while the second one can be interpreted as an
elastic impact of the system against $\mathcal{D}$
(see~\cite{Ibort_de_Leon_Lacomba_Marrero_Martin:Geometric_formulation_of_Carnots_theorem}).
This is what will provide the preservation of energy. Note that we
can combine both equations into
\[D_1 L_d(q_k, q_{k+1})+(\mathcal{P}^*-\mathcal{Q}^*)D_2 L_d (q_{k-1},q_k)=0,\]
from which we see that the system defines a unique discrete
evolution operator if and only if the matrix $(D_{12}L_d)$ is
regular, that is, if the discrete Lagrangian is regular. Locally,
the method can be written as
\begin{subequations}
\begin{align}
D_1 L_d(q_k, q_{k+1})+ D_2 L_d (q_{k-1},
q_k)&=(\lambda_k)_b\,\mu^b
\label{www1}\\
g^{ij}(q_k) \mu^{a}_i(q_k)\left( \frac{\partial L_d}{\partial
q^j_0}(q_k, q_{k+1})- \frac{\partial L_d}{\partial q^j_1}(q_{k-1},
q_k)\right)&=0.\label{www2}
\end{align}
\end{subequations}

Using the discrete Legendre transformations defined above,
define the pre- and post-momenta, which are covectors at $q_k$, by
\begin{align*}
p^+_{k-1, k}&=p^+(q_{k-1}, q_k)= \F^+ L_d(q_{k-1}, q_k)=D_2L_d(q_{k-1}, q_k)\\
p^-_{k, k+1}&=p^-(q_k, q_{k+1})= \F^- L_d(q_k, q_{k+1})=-D_1L_d(q_k, q_{k+1}).
\end{align*}
In these terms, equation \eqref{www2} can be
rewritten as
\[
g^{ij}(q_k) \mu^{a}_i(q_k)\left(\frac{(p^-_{k,
k+1})_j+ (p^+_{k-1, k})_j}{2}\right)=0
\]
which means that the average of post- and
pre-momenta satisfies the constraints. In this
sense the proposed numerical method also
preserves the nonholonomic constraints.

We may rewrite the discrete nonholonomic
equations as
\begin{equation}\label{eq:jump of momenta}
p^-_{k, k+1}=\left({\mathcal P}- {\mathcal
Q}\right)^*_{q_k}(p^+_{k-1,k}).
\end{equation}
We interpret this equation as a jump of momenta
during the nonholonomic evolution. Compare this with the condition $p^-_{k, k+1}=p^+_{k-1, k}$ imposed by the discrete Euler--Lagrange equations (that is, for unconstrained systems). In our method, the momenta are related by a reflection with respect to the image of the projector $\mathcal{P}^*\colon T^*Q\to (\mathcal{D}^\perp)^o$. 
This is illustrated, in the context of Section~\ref{sec:left-invariant lagrangians on lie groups}, in figure~\ref{fig:momenta}.

\subsection{Left-invariant discrete Lagrangians on Lie groups}\label{sec:left-invariant lagrangians on lie groups}
Consider a discrete nonholonomic Lagrangian system on a Lie group
$G$, with a discrete Lagrangian $L_d\colon G\times G\to \R$ that
is invariant with respect to the left diagonal action of $G$ on
$G\times G$
(see~\cite{Bobenko-Suris:Discrete_time_Lagrangian_mechanics_on_Lie_groups_with_an_application_to_the_Lagrange_top,
Marsden_Pekarsky_Shkoller:Discrete_Euler_Poincare_and_Lie_Poisson_equations}).
We do not impose yet any invariance conditions on the distribution
$\mathcal{D}$. If we write $W_k=g_k^{-1}g_{k+1}$, then
we can define the reduced discrete Lagrangian $l_d\colon G\to \R$ as $l_d(W_k)=L_d(g_k,g_{k+1})$. Note that
$Dl_d(W_k)\in T^*_{W_k}G$.

Computing the derivative, we obtain
\[ p^-_{k,k+1}=-D_1L_d(g_k,g_{k+1})=L^*_{g_k^{-1}}R^*_{W_k}Dl_d(W_k), \]
where $L^*$ and $R^*$ are the mappings on $T^*G$ induced by left and right multiplication on the group, respectively (this should not be confused with the Lagrangian $L$). We use this to write
\[
p^+_{k,k+1}=D_2L_d(g_{k},g_{k+1})=L^*_{g_{k}^{-1}}Dl_d(W_{k})=L^*_{g_{k}^{-1}}R^*_{W_{k}^{-1}}L^*_{g_{k}}p^-_{k,k+1}
=R^*_{W_{k}^{-1}}p^-_{k,k+1}.
\]

\begin{figure}
\includegraphics{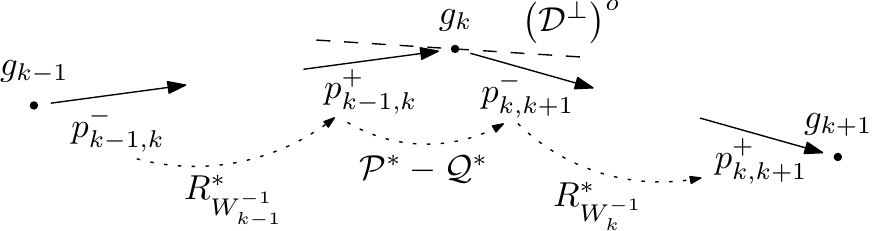}
\caption{Evolution of momenta, depicted here as solid arrows. The right translations are a consequence of the left-invariance of $L_d$, and the reflection at $g_k$ is the proposed method.}
\label{fig:momenta}
\end{figure}

Therefore, the discrete nonholonomic equations \eqref{eq:jump of momenta} become
\begin{equation}
\label{eq:evoluc_en_fibr_cotg}
p^-_{k,k+1}=(\mathcal{P}-\mathcal{Q})^*\left( R^*_{W_{k-1}^{-1}}p^-_{k-1,k}\right).
\end{equation}
The relationships between the pre- and post-momenta are depicted in figure~\ref{fig:momenta}.

Note that we do not need here that the metric used to build the projectors is the metric giving the kinetic energy in the Lagrangian.

\subsection{Left-invariant Lagrangian and projectors}
Take a left-invariant discrete Lagrangian $L_d\colon G\times G\to
\R$ as in the previous section, and assume that $\mathcal{D}$ and $\mathcal{D}^\perp$ are left-invariant. This
is typically a consequence of $\mathcal{D}$ and
the metric on $G$ being left-invariant, although it can be assured
by weaker conditions on the metric (preserving the orthogonality
of $\mathcal{D}_e$ and $\mathcal{D}_e^\perp$ by left
translations). This is equivalent
to the left-invariance of the projectors $\mathcal{P}$ and
$\mathcal{Q}$, which in turn is equivalent to the left-invariance
of $\mathcal{P}-\mathcal{Q}$, as a straightforward verification
shows.

Since our goal is to rewrite equation \eqref{eq:evoluc_en_fibr_cotg} on the dual $\G^*$ of the Lie algebra, we define the discrete body momentum $p_k\colon G\times G \to \G^*$ as
\[ p_k=L^*_{g_k}p^-_{k,k+1}, \]
which agrees with the definition in~\cite{Fedorov_Zenkov:Discrete_nonholonomic_LL_systems_on_Lie_groups}.
Then \eqref{eq:evoluc_en_fibr_cotg} reads
\[ L^*_{g_k^{-1}}p_k=(\mathcal{P}-\mathcal{Q})^*\left( R^*_{W_{k-1}^{-1}}
L^*_{g_{k-1}^{-1}}p_{k-1}
\right). \]
Since $(\mathcal{P}-\mathcal{Q})^*$ is left-invariant, we obtain
\[
p_k=(\mathcal{P}-\mathcal{Q})^*\left( L^*_{g_k}R^*_{W_{k-1}^{-1}}
L^*_{g_{k-1}^{-1}}p_{k-1}\right)
=(\mathcal{P}-\mathcal{Q})^*\left( L^*_{g_{k-1}^{-1}g_k}
R^*_{W_{k-1}^{-1}}
p_{k-1}\right),
\]
that is,
\[ p_k=(\mathcal{P}-\mathcal{Q})^*\left( \Ad^*_{W_{k-1}}
p_{k-1}\right).
\]

\subsection{Preserving energy on Lie groups}
Let us now consider the case where $Q$ is a Lie group $G$, the nonholonomic distribution $\mathcal{D}$ is not necessarily $G$-invariant, and $L$ is regular and bi-invariant.

Since we are restricting ourselves to Lagrangians of mechanical
type, the potential energy is necessarily zero. The
left-invariance of $L$ implies that it must be of the form
\begin{equation}\label{eq:L def por I}
L(v_g)=\frac{1}{2}\left\langle \mathbb{I}g^{-1}v_g,g^{-1}v_g \right\rangle,
\end{equation}
where $\mathbb{I}\colon \G\to \G^*$ is a symmetric non-singular inertia tensor\footnote{In the context of Lie groups, $g$ will denote an element of $G$ instead of the metric.}.
The bi-invariance, however, imposes the equivariance condition $\Ad^*_{g^{-1}}{} \circ \mathbb{I}=\mathbb{I}\circ\Ad_g$ for all $g\in G$, as is straightforward to check. We remark that in this section, the metric used to build the projectors will be the same that defines the Lagrangian. If we take
a discretization $L_d\colon G\times G\to \R$ (which needs to be \textit{left}-invariant only), the equations of motion \eqref{eq:evoluc_en_fibr_cotg} hold. Then we can prove the following result.

\begin{theorem}\label{thm:preservation of energy}
Consider a nonholonomic system on a Lie group with a regular, bi-invariant Lagrangian and with an arbitrary distribution $\mathcal{D}$, and take a discrete Lagrangian that is left-invariant. Then the proposed discrete nonholonomic method~\eqref{eq:propuesta original} is energy-preserving.
\end{theorem}
\begin{proof}
The equivariant inertia tensor $\mathbb{I}$ induces an $\Ad$-invariant scalar product on $\G$ and a bi-invariant metric on $G$. It also defines an inner product $\langle\, ,\rangle_{\mathbb{I}}$ and a corresponding norm $\| \cdot  \|_{\mathbb{I}}$ on each fiber of $T^*G$ that inherit this bi-invariance. If $p\mapsto p^\sharp$ is the index-raising operation associated to the kinetic energy metric, then
\[ \| p_g  \|_{\mathbb{I}}^2=\left\langle p_g,p_g^{\sharp} \right\rangle=\left\langle p_g,L_g\mathbb{I}^{-1}L^*_g p_g \right\rangle=\left\langle p_g,R_g\mathbb{I}^{-1}R^*_g p_g \right\rangle.
 \]
The dual applications of the projectors $\mathcal{P}$ and $\mathcal{Q}$ are orthogonal complementary projectors with respect to this inner product, and thus for $p\in T^*G$,
\[ \| (\mathcal{P}-\mathcal{Q})^*p\|^2_{\mathbb{I}}=
\langle \mathcal{P}^*p,\mathcal{P}^*p \rangle_{\mathbb{I}}+
\langle \mathcal{Q}^*p,\mathcal{Q}^*p \rangle_{\mathbb{I}}=
\| (\mathcal{P}+\mathcal{Q})^*p\|^2_{\mathbb{I}}=
\| p \|^2_{\mathbb{I}}.
\]

The energy function is given in the continuous setting by $H=\left\langle \partial L/ \partial\dot g,\dot g \right\rangle-L$ as a function of the position $g$ and momentum $p=\partial L/ \partial\dot g$. For $L$ given by~\eqref{eq:L def por I} we have
\[ H(g,p)=\frac{1}{2}\left\langle L^*_gp,\mathbb{I}^{-1}L^*_gp \right\rangle=\frac{1}{2}\| p \|^2_{\mathbb{I}}. \]
Proving that the energy is preserved amounts to showing that equation~\eqref{eq:evoluc_en_fibr_cotg} preserves $\| \cdot  \|_{\mathbb{I}}$. Since $\| \cdot  \|_{\mathbb{I}}$ is in particular right-invariant, then $R^*_{W_{k-1}^{-1}}\colon T_{g_{k-1}}^*G\to T_{g_k}^*G$ is an isometry. In addition, we have shown above that $(\mathcal{P}-\mathcal{Q})^*$ is also norm-preserving, so we obtain
\[ H(g_k,p^-_{k,k+1})=H(g_{k-1},p^-_{k-1,k}).\qedhere \]
\end{proof}

\begin{remark}
While the proof above shows that the norm of the post-momenta is preserved, the norm of the pre-momenta is also preserved since they are related by a reflection (equation~\eqref{eq:jump of momenta}).
\end{remark}

\subsection{The average momentum}
Take a discrete nonholonomic system on $G$ as in the previous section, but add the condition that $\mathcal{D}$ is right-invariant. Since the metric on the group is right-invariant, so is the projector $\mathcal{P}$. Take a trajectory of the system and define at each $g_k$ the average momentum
\begin{equation}\label{eq:average momentum}
\widetilde p_k=\frac{1}{2}\left( p^+_{k-1,k}+p^-_{k,k+1} \right).
\end{equation}

Using \eqref{eq:jump of momenta}, \eqref{eq:evoluc_en_fibr_cotg} and the fact that $(\mathcal{P}-\mathcal{Q})^*$ is its own inverse, we have
\begin{equation*}\begin{split}
\widetilde p_k&=\frac{1}{2}\left( (\mathcal{P}-\mathcal{Q})^*(p^-_{k,k+1})+p^-_{k,k+1} \right)=\mathcal{P}^*(p^-_{k,k+1})
=\mathcal{P}^*(R^*_{W_{k-1}^{-1}}p^-_{k-1,k})\\
&=R^*_{W_{k-1}^{-1}}\mathcal{P}^*(p^-_{k-1,k})=R^*_{W_{k-1}^{-1}}\widetilde p_{k-1}.
\end{split}\end{equation*}
Since the norm $\| \cdot  \|_{\mathbb{I}}$ on each fiber of $T^*G$ defined in the proof of Theorem~\ref{thm:preservation of energy} is right-invariant, we obtain $\| \widetilde p_k  \|_{\mathbb{I}}=\| \widetilde p_{k-1} \|_{\mathbb{I}}$, so
\[ H(g_k,\widetilde p_{k})=H(g_{k-1},\widetilde p_{k-1}). \]
In addition, by equation~\eqref{eq:jump of momenta}, we have that
$\mathcal{Q}^*(\widetilde p_k)=0$, so $\widetilde p_k$ satisfies
the constraints.

\subsection{Preservation of the nonholonomic momentum map}
Let us recall some concepts regarding symmetries of nonholonomic systems. Suppose that a Lie group $G$ acts on the configuration manifold $Q$.
Define, for each $q\in Q$, the vector subspace $\G^q$ consisting of those elements of $\G$ whose infinitesimal generators at $q$ satisfy the nonholonomic constraints, i.e.,
\[ \G^q=\left\{ \xi\in \G \,|\, \xi_Q(q)\in \mathcal{D}_q \right\}. \]
The (generalized) bundle over $Q$ whose fiber at $q$ is $\G^q$ is denoted by $\G^\mathcal{D}$.

A horizontal symmetry is an element $\xi\in \G$ such that $\xi_Q(q)\in \mathcal{D}_q$ for all $q\in Q$. Note that a horizontal symmetry is related naturally to a constant section of $\G^\mathcal{D}$.

Now consider a discrete Lagrangian $L_d\colon Q\times Q\to \R$, and define the discrete nonholonomic momentum map $J^{\mathrm{nh}}_d\colon Q\times Q\to (\G^\mathcal{D})^*$ as in~\cite{Cortes_Martinez:Non-holonomic_integrators} by
\begin{align*}
J^{\mathrm{nh}}_d(q_{k-1},q_k)\colon \G^{q_k}&\to \R\\
\xi&\mapsto \left\langle D_2L_d(q_{k-1},q_k),\xi_Q(q_k) \right\rangle.
\end{align*}
For any smooth section $\widetilde \xi$ of $\G^\mathcal{D}$ we have a function $(J^{\mathrm{nh}}_d)_{\widetilde \xi}\colon Q\times Q\to \R$, defined as $(J^{\mathrm{nh}}_d)_{\widetilde \xi}(q_{k-1},q_k)=J^{\mathrm{nh}}_d(q_{k-1},q_k)\left( \widetilde \xi(q_k) \right)$. We can now prove the following result.

\begin{theorem}\label{thm:momentum equation}
Assume that $L_d$ is $G$-invariant, and let $\widetilde \xi$ be a smooth section of $\G^\mathcal{D}$. Then, under the proposed nonholonomic integrator,  $(J^{\mathrm{nh}}_d)_{\widetilde \xi}$ evolves according to the equation
\begin{equation*}
(J^{\mathrm{nh}}_d)_{\widetilde \xi}(q_k,q_{k+1})-(J^{\mathrm{nh}}_d)_{\widetilde \xi}(q_{k-1},q_k)=
\left\langle D_2L_d(q_k,q_{k+1}), \left( \xi_{k+1}-\xi_k \right)_Q(q_{k+1}) \right\rangle
\end{equation*}
where $\xi_k,\xi_{k+1}\in \G$ are the result of dropping the base points of $\widetilde \xi(q_k)$ and $\widetilde \xi(q_{k+1})$ respectively.
\end{theorem}
\begin{proof}
By the invariance of $L_d$ we have
\[ L_d(\exp(s\xi_k) q_k,\exp(s\xi_k) q_{k+1})=L_d(q_k,q_{k+1}), \]
and differentiating at $s=0$ we get
\[\left\langle D_1L_d(q_k,q_{k+1}),(\xi_k)_Q(q_k) \right\rangle +\left\langle D_2L_d(q_k,q_{k+1}),(\xi_k)_Q(q_{k+1}) \right\rangle=0. \]
On the other hand, the proposed integrator implies
\[
(\mathcal P - \mathcal{Q})^*(D_1 L_d(q_k, q_{k+1}))+
D_2 L_d (q_{k-1},
q_k)=0.
\]
{}From this, and using the fact that $(\xi_k)_Q(q_k)\in \mathcal{D}$, we have
\begin{equation*}\begin{split}
(J^{\mathrm{nh}}_d)_{\widetilde \xi}&(q_{k-1},q_k)=\left\langle D_2L_d(q_{k-1},q_k),(\xi_k)_Q(q_k) \right\rangle\\
&=-\left\langle D_1L_d(q_k,q_{k+1}),(\mathcal P - \mathcal{Q})\left( (\xi_k)_Q(q_k) \right) \right\rangle
=-\left\langle D_1L_d(q_k,q_{k+1}), (\xi_k)_Q(q_k)  \right\rangle\\
&=\left\langle D_2L_d(q_k,q_{k+1}),(\xi_k)_Q(q_{k+1})\right\rangle.
\end{split}\end{equation*}
Then
\begin{equation*}\begin{split}
(J^{\mathrm{nh}}_d)_{\widetilde \xi}&(q_k,q_{k+1})-(J^{\mathrm{nh}}_d)_{\widetilde \xi}(q_{k-1},q_k)=\\
&=\left\langle D_2L_d(q_k,q_{k+1}),(\xi_{k+1})_Q(q_{k+1}) \right\rangle-
\left\langle D_2L_d(q_k,q_{k+1}),(\xi_k)_Q(q_{k+1})\right\rangle\\
&=\left\langle D_2L_d(q_k,q_{k+1}),(\xi_{k+1}-\xi_k)_Q(q_{k+1}) \right\rangle\qedhere.
\end{split}\end{equation*}
\end{proof}

\begin{corollary}\label{cor:preservation of momentum}
If $L_d$ is $G$-invariant and $\xi$ is a horizontal symmetry, then the proposed nonholonomic integrator preserves $(J^{\mathrm{nh}}_d)_\xi$.
\end{corollary}

\section{A theoretical example: nonholonomic version of the St\"{o}rmer--Verlet method}\label{section5}

Consider a continuous nonholonomic system determined by the mechanical Lagrangian $L\colon \R^{2n}\to \R$:
\[
L({q}, \dot{q})=\frac{1}{2} \dot{q}^T M \dot{q} -V(q)
\]
(with $M$ a  constant, invertible matrix) and the constraints determined by $\mu(q)\dot{q}=0$ where $\mu(q)$ is a $m\times n$ matrix with $\hbox{rank } \mu=m$.

Consider now the symmetric discretization
\begin{align*}
L_d(q_k, q_{k+1}) &= \frac{1}{2}h L\left(q_k, \frac{q_{k+1}-q_k}{h}\right)+\frac{1}{2}h L\left(q_{k+1}, \frac{q_{k+1}-q_k}{h}\right)\\
                  &= \frac{1}{2h}\left(q_{k+1}-q_k\right)^T M\left(q_{k+1}-q_k\right)-\frac{h}{2}\left(V(q_k)+V(q_{k+1})\right)\,.
\end{align*}
After some straightforward computations we obtain that equations (\ref{www1}) and (\ref{www2}) for the proposed nonholonomic discrete system are
\begin{subequations}\label{poiu}
\begin{align}
q_{k+1}-2q_k+q_{k-1}&=-h^2 M^{-1} \left(V_q(q_{k})+\mu^T(q_k)\widetilde\lambda_k\right)\label{zxc1}\\
0&=\mu(q_k)\left(\frac{q_{k+1}-q_{k-1}}{2h}\right)\label{zxc2},
\end{align}
\end{subequations}
where $V_q(q)=(\partial V/\partial q^i(q))$ and the Lagrange multipliers relate to those in equation (\ref{www1}) by $\widetilde \lambda_k=\lambda_k/h$.
We recognize this set of equations as an obvious extension of the SHAKE method proposed by \cite{Ryckaert} to the case of nonholonomic constraints. The SHAKE method is a generalization of the classical St\"{o}rmer--Verlet method in presence of holonomic constraints. Equations (\ref{poiu}) were proposed by R.\ McLachlan and M.\ Perlmutter \cite{McLachlan_Perlmutter:Integrators_for_nonholonomic_mechanical_systems} (see equations (5.3) therein) as a reversible method for nonholonomic systems \emph{not based} in the Discrete Lagrange--d'Alembert principle.

The momentum components are approximated by the average momentum
$\widetilde{p}_k=M(q_{k+1}-q_{k-1})/2h$ given by equation~\eqref{eq:average momentum}. Denoting $p_{k+1/2}=M(q_{k+1}-q_k)/h$,  equations (\ref{zxc1}) and (\ref{zxc2}) are now rewritten in the form
\begin{align*}
p_{k+1/2}&=\widetilde{p}_k -\frac{h}{2}\left(  V_q(q_{k})+\mu^T(q_k)\widetilde{\lambda}_k\right),\\
q_{k+1}&= q_k+hM^{-1}p_{k+1/2},\\  
0&=\mu(q_k)M^{-1} \widetilde{p}_k. 
\end{align*}

The definition of  $\widetilde{p}_{k+1}$  requires the knowledge of $q_{k+2}$ and, therefore, it  is is natural to apply another step of the algorithm
(\ref{www1}) and (\ref{www2}) to avoid this difficulty. Then, we obtain the new equations:
\begin{align*}
\widetilde{p}_{k+1}&={p}_{k+1/2} -\frac{h}{2}\left( V_q(q_{k+1})+\mu^T(q_{k+1})\widetilde{\lambda}_{k+1}\right),\\
0&=\mu(q_{k+1})M^{-1} \widetilde{p}_{k+1}.
\end{align*}

The interesting result is that we obtain a natural extension of the RATTLE algorithm for holonomic systems to the case of nonholonomic systems. Unifying the equations above we obtain the following numerical scheme
\begin{subequations}\label{wer}
\begin{align}
p_{k+1/2}&=\widetilde{p}_k -\frac{h}{2}\left( V_q(q_{k})+\mu^T(q_k)\widetilde{\lambda}_k\right),\label{azxc1111}\\
q_{k+1}&= q_k+hM^{-1}p_{k+1/2},  \label{azxc1211}\\
0&=\mu(q_k)M^{-1} \widetilde{p}_k\label{azxc2211},\\
\widetilde{p}_{k+1}&={p}_{k+1/2} -\frac{h}{2}\left(V_q(q_{k+1})+\mu^T(q_{k+1})\widetilde{\lambda}_{k+1}\right),\label{azxc11111}\\
0&=\mu(q_{k+1})M^{-1} \widetilde{p}_{k+1}\label{azxc22211}.
\end{align}
\end{subequations}
These equations allow us to take a triple  $(q_k, \widetilde p_k, \widetilde{\lambda}_k)$ satisfying the constraint equations (\ref{azxc2211}), compute $p_{k+1/2}$ using \eqref{azxc1111} and then $q_{k+1}$ using \eqref{azxc1211}. Then, equations \eqref{azxc11111} and \eqref{azxc22211} are used to compute the remaining components of the triple $(q_{k+1}, \widetilde p_{k+1}, \widetilde{\lambda}_{k+1})$.
Of course, from Theorem \ref{thm:preservation of energy} we obtain that, in the case $V=0$, the numerical method is energy preserving.

\begin{remark}
{}From this Hamiltonian point of view, we have shown that the initial
conditions for  this numerical scheme are constrained in a natural way ($(q_0, \widetilde{p}_0)$ with  $\mu(q_0)M^{-1} \widetilde{p}_0=0$), that is, the initial conditions are exactly the same as those for  the continuous system.
However, if we want to maintain the algorithm in the cartesian product $Q\times Q$, then the appropriate set of initial conditions is now
\begin{align}
{\mathcal M}_0&=\{ (q_0, q_1)\in Q\times Q\; |\;  \F^- L_d(q_0, q_{1})\in(\mathcal{D}^\perp)^o\} \label{initial-cond}\\
&=\left\{ (q_0, q_1)\in Q\times Q\; \middle|\; g^{ij}(q_0) \mu^{a}_i(q_0)\frac{\partial L_d}{\partial q^j_0}(q_0, q_{1})=0\right\}. \nonumber
\end{align}
In the particular case of the nonholonomic projection of the St\"{o}rmer--Verlet method we have that
\begin{equation*}
{\mathcal M}_0=\left\{ (q_0, q_1)\in \R^n\times \R^n\; \middle|\;    \mu(q_0) M^{-1}\left( M\frac{q_{1}-q_{0}}{h} +\frac{h}{2}V_q(q_0)\right)=0\right\} .
\end{equation*}
Thus, if $(q_0,q_1)\in \mathcal{M}_0$, we define
\[
\widetilde{p}_0=\frac{\partial L_d}{\partial q_0}(q_0,q_1)=M\frac{q_{1}-q_{0}}{h} +\frac{h}{2}V_q(q_0)=p_{0+1/2}+\frac{h}{2}V_q(q_0)
\]
{}From the expression of $\mathcal{M}_0$ we have that \eqref{azxc2211} holds for $k=0$, and the definition of $p_{0+1/2}$ yields precisely equation \eqref{azxc1211}. If we take $\widetilde \lambda=0$ then \eqref{azxc1111} holds too. Thus, $(q_0, \widetilde p_0, 0)$ can be used to initialize the algorithm \eqref{wer}.
\end{remark}

\begin{remark}
In the particular case where the constraints are integrable, that is, the motion is only defined on a submanifold $N$ of $Q$, then the most natural choice  is  to restrict the discrete Lagrangian to $N\times N$: $(L_d)_{|N\times N}$ (see \cite{Marsden_West:Discrete_mechanics_and_variational_integrators} and references therein).
In a local description $N$ is determined by the vanishing of a family of independent  functions $g^a(q)=0$, $1\leq a\leq m$. Differentiating, we obtain new constraints
\begin{equation}\label{poi}
\frac{\partial g^a}{\partial q^i}(q)\dot{q}^i=0
\end{equation}
 which are satisfied by the trajectories $(c(t), \dot{c}(t))$ in the continuous problem.

If we directly apply our method to a holonomic system we obtain  the preservation of constraints \eqref{poi} but the computed numerical solution will not usually lie on the constraint submanifold $g^a(q)=0$. For instance, it seems more natural to change \eqref{azxc2211} by $g^a(q_{k+1})=0$, as appears in the classical RATTLE method. Nevertheless, in the case $V=0$, our method has as an additional feature the preservation of energy.
We could say that the proposed method is specifically designed for nonintegrable constraints.
\end{remark}

\section{Numerical examples}\label{section6}

\begin{example}

The following typical example will illustrate
some of the constructions of previous sections.
It corresponds to a discretization of the
nonholonomic particle in $\R^3$ described by
\[
L(x, y, z, \dot{x}, \dot{y}, \dot{z})=
\frac{1}{2}\left( \dot{x}^2+
\dot{y}^2+\dot{z}^2\right)
\]
and the nonholonomic constraint
$\varphi=\dot{z}-y\dot{x}=0$, which is represented by the distribution
\[
{\mathcal D}=\text{span}
\left\{\frac{\partial}{\partial
x}+y\frac{\partial}{\partial z},
\frac{\partial}{\partial y}\right\}.
\]
Lagrange--d'Alembert's principle gives the equations of motion
\begin{align*}
\ddot{x}+y\ddot{z}&=0\\
\ddot{y}&=0\\
\dot{z}-y\dot{x}&=0.
\end{align*}

Discretize the system by defining the discrete Lagrangian $L_d\colon\R^3\times \R^3\to \R$ as
\[
L_d(x_0, y_0, z_0, x_1, y_1,
z_1)=\frac{1}{2}\left[
\left(\frac{x_1-x_0}{h}\right)^2+
\left(\frac{y_1-y_0}{h}\right)^2 +
\left(\frac{z_1-z_0}{h}\right)^2\right].
\]
Then the discrete nonholonomic equations are
\begin{subequations}
\begin{align}
\left( \frac{x_2-2x_1+x_0}{h^2} \right) +
y_1\left(
\frac{z_2-2z_1+z_0}{h^2} \right)&=0\label{eq:particula:1}\\
\frac{y_2-2y_1+y_0}{h^2}&=0\label{eq:particula:2}\\
\frac{z_2-z_0}{2h}-y_1\frac{x_2-x_0}{2h}&=0.\label{eq:particula:3}
\end{align}
\end{subequations}
Regarding $\R^3$ as a Lie group under translations, the Euclidean metric is bi-invariant. Since $L$ is induced by this metric and $L_d$
is left-invariant, we have preservation of energy by Theorem~\ref{thm:preservation of energy}. Figure~\ref{fig:energy_behaviour_particle} compares the energy behavior for our method against the DLA (discrete Lagrange--d'Alembert) algorithm in~\cite{Cortes_Martinez:Non-holonomic_integrators}.

In order to write the discrete nonholonomic momentum equation in Theorem~\ref{thm:momentum equation} with respect to this group action, take two linearly independent sections of $\G^\mathcal{D}$ given by $\widetilde \xi_1(x,y,z)=(1,0,y)$ and $\widetilde \xi_2(x,y,z)=(0,1,0)$. The equation for $\widetilde \xi_1$ reads
\[
\left(\frac{x_2-x_1}{h^2}+y_2
\frac{z_2-z_1}{h^2}\right)-\left(\frac{x_1-x_0}{h^2}+y_1
\frac{z_1-z_0}{h^2}\right)=(y_2-y_1)\left(\frac{z_2-z_1}{h^2}\right),
\]
which turns out to be \eqref{eq:particula:1}.
Similarly, if we consider $\widetilde \xi_2$ we reobtain~\eqref{eq:particula:2}.

The DLA method proposed in~\cite{Cortes_Martinez:Non-holonomic_integrators} also yields equations~\eqref{eq:particula:1} and~\eqref{eq:particula:2},
which is reasonable since both methods fulfill the discrete nonholonomic momentum equation. However, the DLA method replaces~\eqref{eq:particula:3} by a discretization of the constraints that does not involve $(x_0,y_0,z_0)$, such as
\[ \frac{z_2-z_1}{h}-\left(\frac{y_2+y_1}{2}\right)\frac{x_2-x_1}{h}=0. \]
\begin{figure}
\includegraphics[scale=.7]{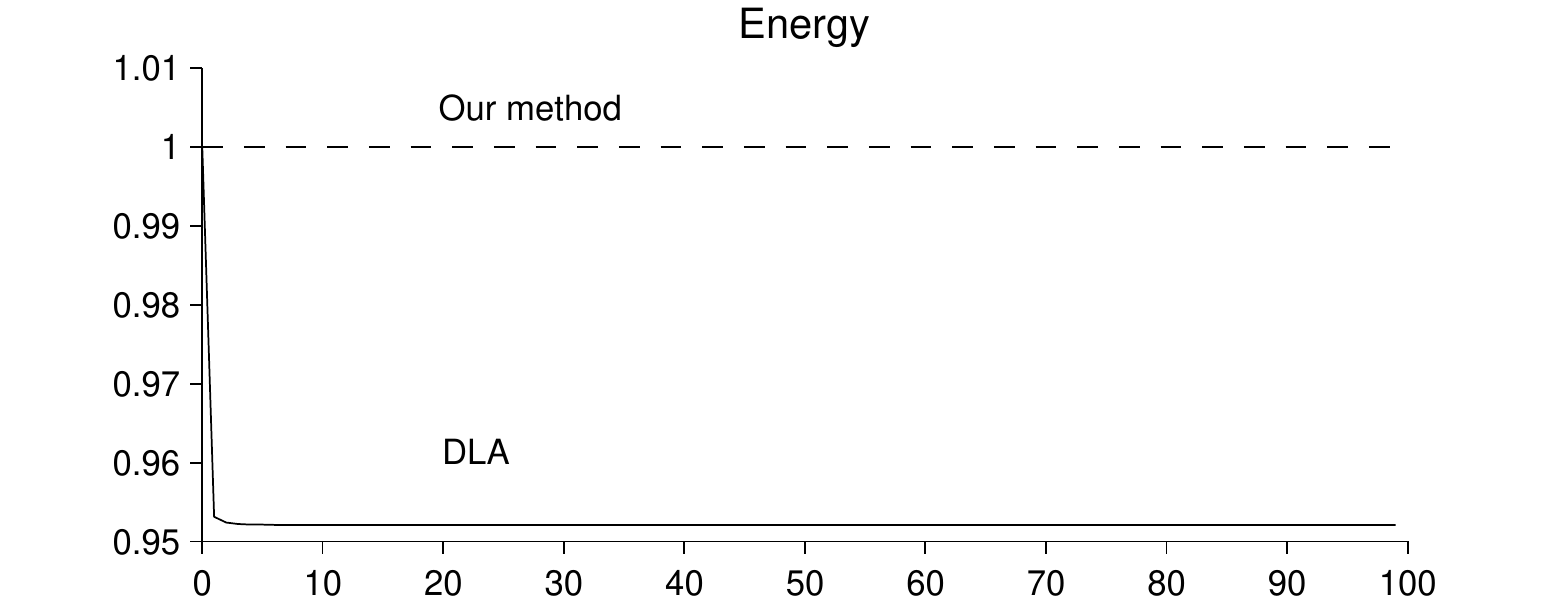}
\caption{Energy behaviour for the nonholonomic particle using our method and the DLA method in~\cite{Cortes_Martinez:Non-holonomic_integrators}.}
\label{fig:energy_behaviour_particle}
\end{figure}
\end{example}

\begin{example}
The snakeboard is a modified version of the traditional
skateboard, where the rider uses his own momentum, coupled with
the constraints, to move the system. The configuration manifold is
$Q=\mathrm{SE}(2)\times \T^2$ with coordinates $(x, y, \theta, \psi, \phi)$
as in figure~\ref{fig:snakeboard}. The center of the board, which is also the center of mass, is located at $(x,y)$. We are considering here the
case where the angles of the front and rear wheel axles are equal
and opposite, as
in~\cite{Bullo_Lewis:Geometric_control_of_mechanical_systems,Lewis:Simple_mechanical_control_systems_with_constraints}. However, we measure these angles with respect to the board instead of the $x$-axis. Figure~\ref{fig:snakeboard} shows a configuration with all the angles positive.

\begin{figure}[hb]
\includegraphics{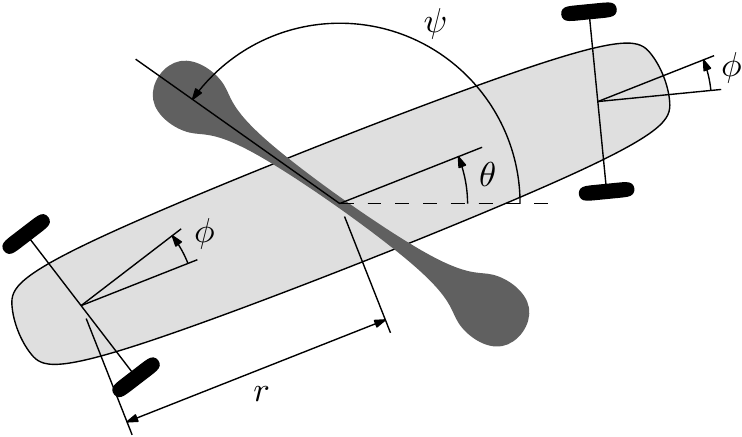}
\caption{The snakeboard. The dashed line is aligned with the $x$-axis (not depicted).}
\label{fig:snakeboard}
\end{figure}

The continuous system is described by the Lagrangian
\begin{eqnarray*}
 L(q, \dot{q})= \frac{1}{2}m (\dot{x}^2+\dot{y}^2) + \frac{1}{2}
 (J+2J_1)\dot{\theta}^2 +\frac{1}{2}J_0\dot{\psi}^2
 +J_1\dot{\phi}^2
 \end{eqnarray*}
where $m$ is the total mass of the system, $J$ is the moment of
inertia of the board about its center,
$J_0$ is the moment of inertia of the rotor
mounted on the board and $J_1$  is the moment of
inertia of each wheel axle about its center. We assume the moments of inertia of the axles about the center of the board to be included in $J$. The distance between the
center of the board and the wheels is denoted by $r$.

The wheels are not allowed to slide sideways, so the constraints turn out to be
\begin{align*}
\dot{x}\sin (\theta+\phi) -\dot{y}\cos (\theta+\phi)
+r\dot{\theta}\cos(\phi)&=0\\
\dot{x}\sin (\theta-\phi) -\dot{y}\cos (\theta-\phi)
-r\dot{\theta}\cos(\phi)&=0.
\end{align*}
If we define the functions $a=r\cos\theta\cos\phi$,
$b=r\sin\theta\cos\phi$ and $c=-\sin\phi$, then
the constraint distribution is
\[
{\mathcal D}=\operatorname{span}\left\{ \frac{\partial}{\partial \psi},
\frac{\partial}{\partial \phi},
a\frac{\partial}{\partial x}+b\frac{\partial}{\partial y}+c\frac{\partial}{\partial \theta}\right\}.
\]
Endow $Q$ with the Riemannian metric associated to the Lagrangian. This is represented in coordinates by the diagonal matrix
\[ \mathbb{I}=
\operatorname{diag}(m,m,J',J_0,2J_1),
\]
where $J'=J+2J_1$.
The orthogonal complement to $\mathcal{D}$ is then
\[ \mathcal{D}^\perp=\operatorname{span}\left\{ J'c\frac{\partial}{\partial x}-ma\frac{\partial}{\partial \theta},b\frac{\partial}{\partial x}-a\frac{\partial}{\partial y} \right\}.\]
The projection $\mathcal{Q}\colon TQ\to \mathcal{D}^\perp$ is given in coordinates by the matrix
\[ \mathcal{Q}=\frac{1}{J'c^2+m(a^2+b^2)}\begin{bmatrix}
  J'c^2+mb^2 & -mab & -J'ac & 0 & 0 \\
  -mab & J'c^2+ma^2 & -J'bc & 0 & 0 \\
  -mac & -mbc & m(a^2+b^2) & 0 & 0 \\
  0 & 0 & 0 & 0 & 0 \\
  0 & 0 & 0 & 0 & 0 \\
\end{bmatrix}, \]
which depends on $(\theta,\phi)$, and its dual $\mathcal{Q}^*$ is represented by the transpose.

Consider the discretization of this system determined by the discrete Lagrangian
\begin{equation*}\begin{split}
L_d (q_k, q_{k+1})&= \frac{1}{h^2}\left( \frac{1}{2}m (\Delta{x_k^2}+\Delta{y_k^2}) +
\frac{1}{2}
 (J+2J_1)\Delta{\theta_k^2} +\frac{1}{2}J_0\Delta{\psi_k^2}
 +J_1\Delta{\phi_k^2}\right)\\
 &=\frac{1}{2h^2}\Delta q_k^T \mathbb{I} \Delta q_k
\end{split}\end{equation*}
where $q_k=(x_k, y_k, \theta_k, \psi_k, \phi_k)$ (a column vector) and $\Delta
z_k=z_{k+1}-z_k$.

The discrete nonholonomic equations~\eqref{eq:propuesta original} can be written as
\[
D_1 L_d(q_k, q_{k+1})+
(\text{Id} - 2\mathcal{Q})^*_{q_k} (D_2 L_d (q_{k-1},
q_k))=0,
\]
so in matricial form we get
\begin{equation}\label{eq:ecuacion_matricial}
 \frac{1}{h^2}\left(  \mathbb{I}\Delta q_k+(\text{Id}-2\mathcal{Q}^T_{q_k})(-\mathbb{I}\Delta q_{k-1})\right)=0,
\end{equation}
that is,
\[ q_{k+1}=  (\text{Id}-2\mathbb{I}^{-1}\mathcal{Q}^T_{q_k}\mathbb{I})\Delta q_{k-1}+q_k . \]

Regarding the configuration space $\mathrm{SE}(2)\times \mathbb{T}^2$ as a Lie group, $L$ is left-invariant. However, it cannot be right-invariant,
because there are no bi-invariant metrics in $\mathrm{SE}(2)$. If one changes the group structure for the variables $(x,y,\theta)$ from $\mathrm{SE}(2)$ to $\R^2\times \mathrm{S}^1$, then both the continuous and discrete Lagrangians are bi-invariant. The numerical method itself does not depend on which symmetry group one takes, but considering this last group structure allows us to apply Theorem~\ref{thm:preservation of energy} to show that there is preservation of energy.

On the other hand, we can still use the non-abelian group structure to write the discrete nonholonomic momentum equations, since only the left-invariance of $L_d$ is required. Let us consider the action of the subgroup $\mathrm{SE}(2)$ on $\mathrm{SE}(2)\times \mathbb{T}^2$, and
take the typical basis of $\mathfrak{se}(2)$: $e_1=\left(
\begin{smallmatrix}
0&0&1\\
0&0&0\\
0&0&0
\end{smallmatrix}
\right)$, $e_2=\left(
\begin{smallmatrix}
0&0&0\\
0&0&1\\
0&0&0
\end{smallmatrix}
\right)$ and $e_3=\left(
\begin{smallmatrix}
0&-1&0\\
1&0&0\\
0&0&0
\end{smallmatrix}
\right)$. Consider the section $\tilde{\xi}\colon Q\to \mathfrak{se}(2)$ defined by $\tilde{\xi}(x, y, \theta, \psi,
\phi)=(a(\theta, \phi)+c(\theta, \phi)y)e_1+(b(\theta,
\phi)-c(\theta, \phi)x)e_2+c(\theta, \phi)e_3$, so we have
$\tilde{\xi}_Q=a(\theta, \phi)\frac{\partial}{\partial
x}+b(\theta, \phi)\frac{\partial}{\partial y}+c(\theta,
\phi)\frac{\partial}{\partial \theta}$. Therefore, the discrete
nonholonomic momentum equation in this case is
\begin{equation*}
\begin{split}
(J^{\mathrm{nh}}_d)_{\tilde{\xi}}(q_k, q_{k+1})&-
(J^{\mathrm{nh}}_d)_{\tilde{\xi}}(q_{k-1},
q_k)=\\
&m(a(\theta_{k+1}, \phi_{k+1})-a(\theta_{k},
\phi_{k}))\frac{x_{k+1}-x_{k}}{h^2}\\
+&m(b(\theta_{k+1}, \phi_{k+1})-b(\theta_{k},
\phi_{k}))\frac{y_{k+1}-y_{k}}{h^2}\\
+&(J+2J_1)(c(\theta_{k+1}, \phi_{k+1})-c(\theta_{k},
\phi_{k}))\frac{\theta_{k+1}-\theta_{k}}{h^2}.
\end{split}
\end{equation*}
As an additional application, our method is ready to introduce controlled external forces. For instance  we have  added two controls: one applying equal but opposite torques on the wheel axles, and the other one on the rider. This was done by including appropriate terms on the right-hand side of equation~\eqref{eq:ecuacion_matricial}. The figure below shows a simulation where the snakeboard starts from rest and the controls are sinusoidal, with the same phase and frequency. This achieves the typical ``snake-like'' forward motion of the snakeboard, with increasing speed.

\begin{figure}[h]
\includegraphics[scale=.8]{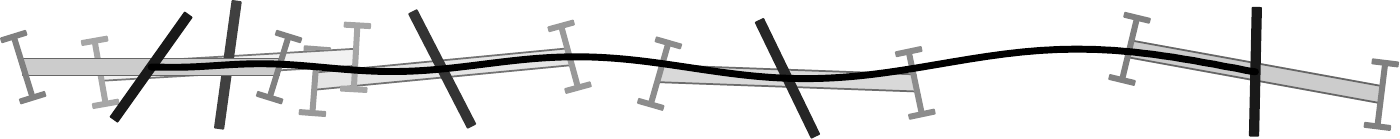}
\caption{The controlled snakeboard, moving left to right.}
\label{fig:controlled_snakeboard}
\end{figure}

\end{example}

\begin{example}
The Chaplygin sleigh consists in a rigid body that moves on a plane and is supported at three points. One of them is a knife edge and cannot slide sideways, and the other two can slide freely. Assume that the sleigh is symmetric, meaning that the center of mass is located on the line determined by the knife edge, at a distance $a$ of the point of contact $(x,y)$ (see figure~\ref{fig:sleigh}).

\begin{figure}[hb]
\includegraphics{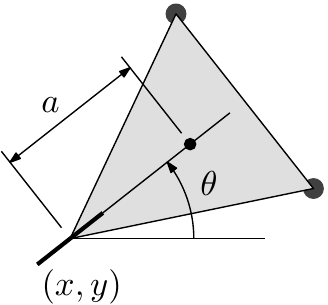}
\caption{The Chaplygin sleigh.}
\label{fig:sleigh}
\end{figure}

The position of the sleigh is determined by $q=(x,y,\theta)\in \R^2\times \mathrm{S}^1$, and the nonholonomic constraint is $\dot x\sin\theta-\dot y\cos\theta=0$. If $m$ is the mass of the sleigh, $I$ is its moment of inertia and $(x_C,y_C)$ denotes the position of the center of mass, then the Lagrangian is
\[ L=\frac{1}{2}m\left( \dot x_C^2+\dot y_C^2 \right)+\frac{1}{2}I\dot \theta^2=\frac{1}{2}m\left(\dot x^2-2a\dot \theta \dot x\sin\theta+\dot y^2+2a \dot \theta\dot y\cos\theta+a^2\dot \theta^2 \right)+\frac{1}{2}I\dot \theta^2. \]
The kinetic energy metric is represented by the matrix
\[ \begin{bmatrix}
  m & 0 & -am\sin\theta \\
  0 & m & am\cos\theta \\
  -am\sin\theta & am\cos\theta & I+ma^2 \\
\end{bmatrix} \]
so the constraint distribution and its orthogonal complement are
\begin{align*}
\mathcal{D}&=\operatorname{span}\left\{ \cos\theta\frac{\partial}{\partial x}+\sin\theta\frac{\partial}{\partial y},\frac{\partial}{\partial \theta} \right\}\\
\mathcal{D}^\perp&=\operatorname{span}\left\{-\sin\theta \frac{\partial}{\partial x}+\cos\theta\frac{\partial}{\partial y}-\frac{am}{I+ma^2}\frac{\partial}{\partial \theta} \right\}.
\end{align*}
The dual of the projector onto $\mathcal{D}^\perp$ is then given by
\[ \mathcal{Q}^*=\begin{bmatrix}
  \sin^2\theta & -\sin\theta\cos\theta & \displaystyle\frac{am\sin\theta}{I+ma^2} \\
  -\sin\theta\cos\theta & \cos^2\theta & \displaystyle-\frac{am\cos\theta}{I+ma^2} \\
  0 & 0 & 0 \\
\end{bmatrix}.\]

Discretize the Lagrangian by replacing $ \dot x$ by $(x_1-x_0)/h$ (analogously for $\dot y$ and $\dot \theta$), and $\theta$ by $(\theta_0+\theta_1)/2$. We have applied the DLA algorithm, discretizing the constraints by $(x_2-x_1)\sin((\theta_1+\theta_2)/2)-(y_2-y_1)\cos((\theta_1+\theta_2)/2)=0$, and compared the results with the trajectory of the continuous system. This trajectory was obtained by applying standard numerical methods to the Lagrange--d'Alembert differential equations (see for example~\cite[p.\ 25]{Bloch:Nonholonomic_mechanics_and_control}).
Figure~\ref{fig:sleigherror} shows the evolution of $( ( x_k-\bar x_k )^2+( y_k-\bar y_k )^2+( \theta_k-\bar \theta_k )^2 )^{1/2}$ for both DLA and our method, where $(\bar x_k,\bar y_k,\bar \theta_k)$ are the values at $t=hk$ of the trajectory of the continuous system. The results shown correspond to a particular trajectory with the initial points extracted from the continuous solution, but in general the errors are similar for the two methods. We used $m=J=1$, $a=.2$, $q_0=(0,0,0)$ and $q_1=(-.2395,-.0070,.0589)$, which produces the heart-shaped loop typically described by the sleigh.

\begin{figure}
\includegraphics[scale=.7]{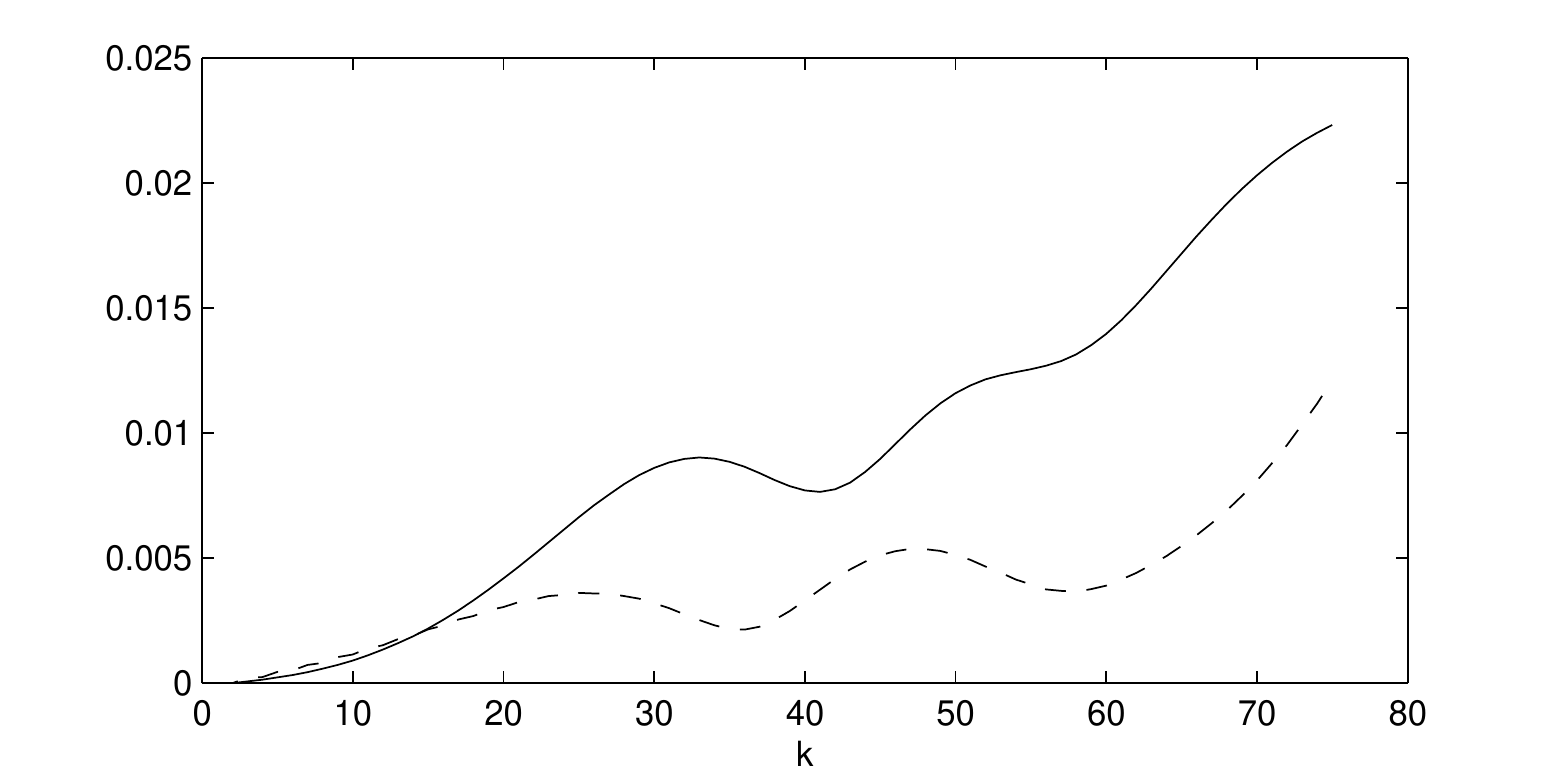}
\caption{Error in $\R^3$ of the trajectories computed with our method (dashed line) and DLA (solid).}
\label{fig:sleigherror}
\end{figure}

It is worth mentioning that if we take a different discretization of the constraints for the DLA algorithm, such as $(x_2-x_1)\sin\theta_1-(y_2-y_1)\cos\theta_1=0$, the error becomes larger by one to two orders of magnitude. Taking the right discretization is crucial in the DLA algorithm; in contrast, the accuracy of our method is close to that of DLA without the need of such a choice.
\end{example}

\section{Conclusions and future work}

In this paper,  we propose a geometric integrator for nonholonomic
mechanical systems for which the constraints are not required to be discretized.
The integrator is different from the usual discrete analogue of the Lagrange--d'Alembert
(DLA) principle which is presented in the works \cite{Cortes_Martinez:Non-holonomic_integrators,McLachlan_Perlmutter:Integrators_for_nonholonomic_mechanical_systems}.
As initial conditions we propose points $(q_0,q_1)$ satisfying \eqref{initial-cond}.

Our method preserves in average the nonholonomic constraints, and the nonholonomic
momentum map is also preserved. In addition, when the configuration space is a Lie
group and some invariance conditions for the continuous and discrete Lagrangians are
satisfied, we prove that the energy is preserved. In the particular case of a typical symmetric
discretization of a mechanical Lagrangian we obtain a natural generalization of the well-known
RATTLE method for holonomic constraints. In addition, several interesting concrete examples
illustrate these results.

Of course,  much work remains to be done to clarify the nature of
discrete nonholonomic mechanics. A large part of this future work was
stated in \cite{McLachlan_Perlmutter:Integrators_for_nonholonomic_mechanical_systems}
and, in particular, we emphasize the following important topics: a complete backward error analysis and
the construction of a discrete exact model for a continuous nonholonomic system;
studying discrete nonholonomic systems that preserve a volume form on the constraint
surface, mimicking the continuous case; analyzing the discrete Hamiltonian framework; and the
construction of integrators depending on different
discretizations.

For the case of reduced systems, it is possible to adapt the Lie-groupoid techniques introduced
in the papers \cite{Iglesias_Marrero_MartindeDiego_Martinez:Discrete_Nonholonomic_Lagrangian_Systems_on_Lie_Groupoids,
Marrero_MartindeDiego_Martinez:Discrete_Lagrangian_and_Hamiltonian_mechanics_on_Lie_groupoids},
considering now a fibred metric on the associated Lie algebroid and the induced orthogonal projectors.

In future works, we will study these problems and, moreover, we will develop explicit constructions
of higher order nonholonomic methods
and applications to numerical methods for optimal control problems (of nonholonomic systems).

\section*{Acknowledgments}
This work has been partially supported by MEC (Spain) Grant  MTM
2007-62478, project ``Ingenio Mathematica'' (i-MATH) No. CSD
2006-00032 (Consolider-Ingenio 2010) and Project SIMUMAT
S-0505/ESP/0158 of the CAM. S.\ Ferraro also wants to thank
SIMUMAT for a Research contract and D.\ Iglesias to CSIC for a
JAE Research Contract.

The authors would like to thank the referees for the interesting
and helpful comments which have helped
to improve the contents of the paper.

\bibliography{allrefs}

\newcommand\oneletter[1]{#1}\newcommand\Yu{Yu}\newcommand{\ct}{t}\def\cprime{$%
'$}
\begin{thebibliography}{10}

\bibitem{Abraham_Marsden:Foundations_of_Mechanics}
Ralph Abraham and Jerrold~E. Marsden.
\newblock {\em Foundations of mechanics}.
\newblock Benjamin/Cummings Publishing Co. Inc. Advanced Book Program, Reading,
  Mass., 1978.
\newblock Second edition, revised and enlarged, with the assistance of Tudor
  Ra{\ct}iu and Richard Cushman.

\bibitem{Bloch:Nonholonomic_mechanics_and_control}
Anthony~M. Bloch.
\newblock {\em Nonholonomic mechanics and control}, volume~24 of {\em
  Interdisciplinary Applied Mathematics}.
\newblock Springer-Verlag, New York, 2003.
\newblock With the collaboration of J. Baillieul, P. Crouch and J. Marsden,
  With scientific input from P. S. Krishnaprasad, R. M. Murray and D. Zenkov,
  Systems and Control.

\bibitem{Bloch_Krishnaprasad_Marsden_Ratiu:The_Euler-Poincare_equations_and_do%
uble_bracket_dissipation}
Anthony~M. Bloch, P.~S. Krishnaprasad, Jerrold~E. Marsden, and Tudor~S.
  Ra{\ct}iu.
\newblock The {E}uler-{P}oincar\'e equations and double bracket dissipation.
\newblock {\em Comm. Math. Phys.}, 175(1):1--42, 1996.

\bibitem{Bloch_Marsden_Zenkov:Nonholonomic_dynamics}
Anthony~M. Bloch, Jerrold~E. Marsden, and Dmitry~V. Zenkov.
\newblock Nonholonomic dynamics.
\newblock {\em Notices Amer. Math. Soc.}, 52(3):324--333, 2005.

\bibitem{Bobenko-Suris:Discrete_time_Lagrangian_mechanics_on_Lie_groups_with_a%
n_application_to_the_Lagrange_top}
Alexander~I. Bobenko and {\oneletter{\Yu}}ri~B. Suris.
\newblock Discrete time {L}agrangian mechanics on {L}ie groups, with an
  application to the {L}agrange top.
\newblock {\em Comm. Math. Phys.}, 204(1):147--188, 1999.

\bibitem{Bullo_Lewis:Geometric_control_of_mechanical_systems}
Francesco Bullo and Andrew~D. Lewis.
\newblock {\em Geometric control of mechanical systems}, volume~49 of {\em
  Texts in Applied Mathematics}.
\newblock Springer-Verlag, New York, 2005.
\newblock Modeling, analysis, and design for simple mechanical control systems.

\bibitem{Cantrijn_Cortes_deLeon_MartindeDiego:On_the_geometry_of_generalized_C%
haplygin_systems}
Frans Cantrijn, Jorge Cort{\'e}s, Manuel de~Le{\'o}n, and David
  Mart{\'{\i}}n~de Diego.
\newblock On the geometry of generalized {C}haplygin systems.
\newblock {\em Math. Proc. Cambridge Philos. Soc.}, 132(2):323--351, 2002.

\bibitem{Cendra_Ibort_de_Leon_de_Diego:A_generalization_of_Chetaevs_principle_%
for_a_class_of_higher_order_nonholonomic_constraints}
Hern{\'a}n Cendra, Alberto Ibort, Manuel de~Le{\'o}n, and David
  Mart{\'{\i}}n~de Diego.
\newblock A generalization of {C}hetaev's principle for a class of higher order
  nonholonomic constraints.
\newblock {\em J. Math. Phys.}, 45(7):2785--2801, 2004.

\bibitem{Cortes:Energy_conserving_nonholonomic_integrators}
Jorge Cort{\'e}s.
\newblock Energy conserving nonholonomic integrators.
\newblock {\em Discrete Contin. Dyn. Syst.}, (suppl.):189--199, 2003.
\newblock Dynamical systems and differential equations (Wilmington, NC, 2002).

\bibitem{Cortes_Martinez:Non-holonomic_integrators}
Jorge Cort{\'e}s and Sonia Mart{\'{\i}}nez.
\newblock Non-holonomic integrators.
\newblock {\em Nonlinearity}, 14(5):1365--1392, 2001.

\bibitem{de_Leon_Martin_de_Diego:On_the_geometry_of_non-holonomic_Lagrangian_s%
ystems}
Manuel de~Le{\'o}n and David Mart{\'\i}n~de Diego.
\newblock On the geometry of non-holonomic {L}agrangian systems.
\newblock {\em J. Math. Phys.}, 37(7):3389--3414, 1996.

\bibitem{de_Leon_Martin_de_Diego_Santamaria-Merino:Geometric_numerical_integra%
tion_of_nonholonomic_systems_and_optimal_control_problems}
Manuel de~Le{\'o}n, David Mart{\'{\i}}n~de Diego, and Aitor
  Santamar{\'{\i}}a-Merino.
\newblock Geometric numerical integration of nonholonomic systems and optimal
  control problems.
\newblock {\em Eur. J. Control}, 10(5):515--521, 2004.

\bibitem{Fedorov_Zenkov:Discrete_nonholonomic_LL_systems_on_Lie_groups}
Yuri~N. Fedorov and Dmitry~V. Zenkov.
\newblock Discrete nonholonomic {LL} systems on {L}ie groups.
\newblock {\em Nonlinearity}, 18(5):2211--2241, 2005.

\bibitem{Hairer_Lubich_Wanner:Geometric_numerical_integration}
Ernst Hairer, Christian Lubich, and Gerhard Wanner.
\newblock {\em Geometric numerical integration}, volume~31 of {\em Springer
  Series in Computational Mathematics}.
\newblock Springer-Verlag, Berlin, second edition, 2006.
\newblock Structure-preserving algorithms for ordinary differential equations.

\bibitem{Ibort_de_Leon_Lacomba_Marrero_Martin:Geometric_formulation_of_Carnots%
_theorem}
Alberto Ibort, Manuel de~Le{\'o}n, Ernesto~A. Lacomba, Juan~C. Marrero, David
  Mart{\'{\i}}n~de Diego, and Paulo Pitanga.
\newblock Geometric formulation of {C}arnot's theorem.
\newblock {\em J. Phys. A}, 34(8):1691--1712, 2001.

\bibitem{Iglesias_Marrero_MartindeDiego_Martinez:Discrete_Nonholonomic_Lagrang%
ian_Systems_on_Lie_Groupoids}
David Iglesias, Juan~C. Marrero, David Mart{\'{\i}}n~de Diego, and Eduardo
  Mart{\'{\i}}nez.
\newblock Discrete nonholonomic {L}agrangian systems on {L}ie groupoids.
\newblock Preprint arXiv:0704.1543v1, to appear in J.\ Nonlinear Sci., 2007.

\bibitem{Lewis:Affine_connections_and_distributions_with_applications_to_nonho%
lonomic_mechanics}
Andrew~D. Lewis.
\newblock Affine connections and distributions with applications to
  nonholonomic mechanics.
\newblock {\em Rep. Math. Phys.}, 42(1-2):135--164, 1998.
\newblock Pacific Institute of Mathematical Sciences Workshop on Nonholonomic
  Constraints in Dynamics (Calgary, AB, 1997).

\bibitem{Lewis:Simple_mechanical_control_systems_with_constraints}
Andrew~D. Lewis.
\newblock Simple mechanical control systems with constraints.
\newblock {\em IEEE Trans. Automat. Control}, 45(8):1420--1436, 2000.
\newblock Mechanics and nonlinear control systems.

\bibitem{Marrero_MartindeDiego_Martinez:Discrete_Lagrangian_and_Hamiltonian_me%
chanics_on_Lie_groupoids}
Juan~C. Marrero, David Mart{\'{\i}}n~de Diego, and Eduardo Mart{\'{\i}}nez.
\newblock Discrete {L}agrangian and {H}amiltonian mechanics on {L}ie groupoids.
\newblock {\em Nonlinearity}, 19(6):1313--1348, 2006.

\bibitem{Marsden_Pekarsky_Shkoller:Discrete_Euler_Poincare_and_Lie_Poisson_equ%
ations}
Jerrold~E. Marsden, Sergey Pekarsky, and Steve Shkoller.
\newblock Discrete {E}uler-{P}oincar\'e and {L}ie-{P}oisson equations.
\newblock {\em Nonlinearity}, 12(6):1647--1662, 1999.

\bibitem{Marsden_West:Discrete_mechanics_and_variational_integrators}
Jerrold~E. Marsden and Matthew West.
\newblock Discrete mechanics and variational integrators.
\newblock {\em Acta Numer.}, 10:357--514, 2001.

\bibitem{McLachlan_Perlmutter:Integrators_for_nonholonomic_mechanical_systems}
R.~McLachlan and M.~Perlmutter.
\newblock Integrators for nonholonomic mechanical systems.
\newblock {\em J. Nonlinear Sci.}, 16(4):283--328, 2006.

\bibitem{McLachlan_Scovel:A_survey_of_open_problems_in_symplectic_integration}
Robert~I. McLachlan and Clint Scovel.
\newblock A survey of open problems in symplectic integration.
\newblock In {\em Integration algorithms and classical mechanics (Toronto, ON,
  1993)}, volume~10 of {\em Fields Inst. Commun.}, pages 151--180. Amer. Math.
  Soc., Providence, RI, 1996.

\bibitem{Neimark_Fufaev:Dynamics_of_Nonholonomic_Systems}
{\oneletter{\Yu}}ri~I. Ne\u{\i}mark and Nikolai~A. Fufaev.
\newblock {\em Dynamics of Nonholonomic Systems}.
\newblock Translations of Mathematical Monographs, Vol. 33. American
  Mathematical Society, Providence, R.I., 1972.

\bibitem{Ryckaert}
Jean-Paul Ryckaert, Giovanni Ciccotti, and Herman J.~C. Berendsen.
\newblock Numerical integration of the cartesian equations of motion of a
  system with constraint: molecular dynamics of $n$-alkanes.
\newblock {\em J. Comput. Physics}, 23:327--341, 1977.

\end{thebibliography}

\end{document}